# THE HARTLEY TRANSFORM
# IN A FINITE FIELD

R. M. Campello de Souza, H. M. de Oliveira, A. N. Kauffman
CODEC - Grupo de Pesquisas em Comunicações
Departamento de Eletrônica e Sistemas - CTG - UFPE
C.P. 7800, 50711 - 970, Recife - PE , Brasil
Phone: +55-81-271-8210 fax:+55-81-271-8215
e-mail: Ricardo@npd.ufpe.br , HMO@npd.ufpe.br,

*Abstract – In this paper, the k-trigonometric functions over the Galois Field GF(q) are introduced and their main properties derived. This leads to the definition of the $cas_k(.)$ function over GF(q), which in turn leads to a finite field Hartley Transform . The main properties of this new discrete transform are presented and areas for possible applications are mentioned.*

## 1. Introduction

Discrete transforms play a very important role in engineering. A significant example is the well known Discrete Fourier Transform (DFT), which has found many applications in several areas, specially in Electrical Engineering. A DFT for finite fields was introduced by Pollard in 1971 [1] and applied as a tool to perform discrete convolutions using integer arithmetic. Since then several new applications of the Finite Field Fourier Transform (FFFT) have been found, not only in the fields of digital signal and image processing [2-5], but also in different contexts such as error control coding and cryptography [6-8].

A second relevant example concerns the Discrete Hartley Transform (DHT) [9], the discrete version of the integral transform introduced by R. V. L. Hartley in [10]. Although seen initially as a tool with applications only on the numerical side and having connections to the physical world only via the Fourier transform, the DHT has proven over the years to be a very useful instrument with many interesting applications [11-13].

In this paper the DHT over a finite field is introduced. In order to obtain a transform that holds some resemblance with the DHT, it is necessary firstly to establish the equivalent of the sinusoidal functions cos and sin over a finite structure. Thus, in the next section, the k-trigonometric functions $\cos_k$ and $\sin_k$ are defined from which the $cas_k$ (**c**osine **a**nd **s**ine) function is obtained and used, in section 3, to introduce a symmetrical discrete transform pair, the finite field Hartley transform, or FFHT for short. A number of properties of the FFHT is presented, including the cyclic convolution property and Parseval's relation. In section 4, the condition for valid spectra, similar to the conjugacy constraints for the Finite Field Fourier Transform, is given. Section 5 contains a few concluding remarks and some possible areas of applications for the ideas introduced in the paper. The FFHT presented here is different from an earlier proposed Hartley Transform in finite fields [14] and appears to be the more natural one.

## 2. k-Trigonometric Functions

The set G(q) of gaussian integers over GF(q) defined below plays an important role in the ideas introduced in this paper (hereafter the symbol := denotes *equal by definition*).

<u>Definition 1:</u> $G(q) := \{a + jb ,\ a, b \in GF(q)\}$, $q = p^r$, r being a positive integer, p being an odd prime for which $j^2 = -1$ is a quadratic non-residue in GF(q), is the set of gaussian integers over GF(q).

Let $\otimes$ denote the cartesian product. It can be shown, as indicated below, that the set G(q) together with the operations $\oplus$ and $*$ defined below, is a field.

<u>Proposition 1</u> : Let

$$\oplus : G(q) \otimes G(q) \rightarrow G(q)$$

$$(a_1 + jb_1 , a_2 + jb_2) \rightarrow (a_1 + jb_1) \oplus (a_2 + jb_2) = (a_1 + a_2) + j(b_1 + b_2)$$

and

$$* : G(q) \otimes G(q) \rightarrow G(q)$$

$(a_1 + jb_1 , a_2 + jb_2) \rightarrow (a_1 + jb_1) * (a_2 + jb_2) =$
$= (a_1a_2 - b_1b_2) + j(a_1b_2 + a_2b_1.$

The structure GI(q) := < G(q) , ⊕ , ∗ > is a field. In fact, GI(q) is isomorphic to $GF(q^2)$. ð

Trigonometric functions over the elements of a Galois field can be defined as follows.

Definition 2 : Let $\alpha$ have multiplicative order N in GF(q), $q = p^r$ , $p \neq 2$. The GI(q)-valued k-trigonometric functions of $\angle(\alpha^i)$ in GF(q) (by analogy, the trigonometric functions of k times the ″angle″ of the ″complex exponential″ $\alpha^i$ ) are defined as

$$\cos_k(\angle\alpha^i) := \frac{1}{2}(\alpha^{ik} + \alpha^{-ik})$$

and

$$\sin_k(\angle\alpha^i) := \frac{1}{2j}(\alpha^{ik} - \alpha^{-ik}),$$

for i , k = 0, 1,..., N-1 .

For simplicity suppose $\alpha$ to be fixed. We write $\cos_k(\angle\alpha^i)$ as $\cos_k(i)$ and $\sin_k(\angle\alpha^i)$ as $\sin_k(i)$. The k-trigonometric functions satisfy properties P1-P8 below. Proofs are straightforward and are omitted here.

P1. Unit Circle: $\sin_k^2(i) + \cos_k^2(i) = 1.$

P2. Even / Odd: $\cos_k(i) = \cos_k(-i)$

$\sin_k(i) = -\sin_k(-i).$

P3. Euler Formula : $\alpha^{ik} = \cos_k(i) + j\sin_k(i).$

P4. Addition of Arcs :

$$\cos_k(i + t) = \cos_k(i)\cos_k(t) - \sin_k(i)\sin_k(t),$$

$$\sin_k(i + t) = \sin_k(i)\cos_k(t) + \sin_k(t)\cos_k(i).$$

P5. Double Arc:

$$\cos_k^2(i) = \frac{1 + \cos_k(2i)}{2}$$

$$\sin_k^2(i) = \frac{1 - \cos_k(2i)}{2}$$

P6. Symmetry :

$$\cos_k(i) = \cos_i(k)$$

$$\sin_k(i) = \sin_i(k).$$

P7. $\cos_k(i)$ Summation:

$$\sum_{k=0}^{N-1} \cos_k(i) = \begin{cases} N , & i = 0 \\ 0 , & i \neq 0 \end{cases}.$$

P8. $\sin_k(i)$ Summation:

$$\sum_{k=0}^{N-1} \sin_k(i) = 0.$$

A simple example is given to illustrate the behavior of such functions.

Example 1 - Let $\alpha = 3$, a primitive element of GF(7). The $\cos_k$ (i) and $\sin_k$ (i) functions take the following values in GF(7):

$\cos_k$ (i)

| (k) \ (i) | 0 | 1 | 2 | 3 | 4 | 5 |
|---|---|---|---|---|---|---|
| 0 | 1 | 1 | 1 | 1 | 1 | 1 |
| 1 | 1 | 4 | 3 | 6 | 3 | 4 |
| 2 | 1 | 3 | 3 | 1 | 3 | 3 |
| 3 | 1 | 6 | 1 | 6 | 1 | 6 |
| 4 | 1 | 3 | 3 | 1 | 3 | 3 |
| 5 | 1 | 4 | 3 | 6 | 3 | 4 |

$\sin_k$ (i)

| (k) \ (i) | 0 | 1 | 2 | 3 | 4 | 5 |
|---|---|---|---|---|---|---|
| 0 | 0 | 0 | 0 | 0 | 0 | 0 |
| 1 | 0 | j | j | 0 | 6j | 6j |
| 2 | 0 | j | 6j | 0 | j | 6j |
| 3 | 0 | 0 | 0 | 0 | 0 | 0 |
| 4 | 0 | 6j | j | 0 | 6j | j |
| 5 | 0 | 6j | 6j | 0 | j | j |

Table 1 – Discrete cosine and sine functions over GF(7).

The k-trigonometric functions have interesting orthogonality properties, such as the one shown in lemma 1.

Lemma 1: The k-trigonometric functions $\cos_k(.)$ and $\sin_k(.)$ are orthogonal in the sense that

$$A := \sum_{k-0}^{N-1} [\cos_k(\angle\alpha^i)\sin_k(\angle\alpha^t)] = 0,$$

where $\alpha$ is an element of multiplicative order N in GF(q).

Proof: By definition 2,

$$A = \sum_{k=0}^{N-1} [\frac{1}{2}(\alpha^{ik} + \alpha^{-ik})\frac{1}{2j}(\alpha^{tk} - \alpha^{-tk})] = = \frac{1}{4j}\sum_{k=0}^{N-1} (\alpha^{k(i+t)} - \alpha^{-k(i+t)} + \alpha^{k(t-i)} - \alpha^{k(i-t)}).$$

Now, If $i = t$, then $A = (0 + 0 + N - N) / 4j = 0$. If $i = -t$, then $A = (N - N + 0 - 0) / 4j = 0$. Otherwise, $A= (0 + 0 + 0 + 0)/4j = 0$. ð

A general orthogonality condition, which leads to a new Hartley Transform, is now presented via the $\text{cas}_k(\angle\alpha^i)$ function. The notation used here follows closely the original one introduced in [10].

Definition 3: Let $\alpha \in$ GF(q), $\alpha \neq 0$. Then

$$\text{cas}_k(\angle\alpha^i) := \cos_k(\angle\alpha^i) + \sin_k(\angle\alpha^i).$$

The set $\{cas_k(.)\}_{k=0, 1,...., N-1}$ , can be viewed as a set of sequences that satisfy the following orthogonality property:

Theorem 1:

$$H := \sum_{k=0}^{N-1} cas_k(\angle\alpha^i)\, cas_k(\angle\alpha^t) = \begin{cases} N , & i = t \\ 0 , & i \neq t \end{cases} ,$$

where $\alpha$ has multiplicative order N.

Proof: From definition 3 it follows that

$$H = \sum_{k=0}^{N-1} [cos_k( i )cos_k( t ) + sin_k( i )sin_k( t ) + +sin_k( i )cos_k( t ) + sin_k( t )cos_k( i )],$$

which, by lemma 1, is the same as

$$H = \sum_{k=0}^{N-1} cos_k( i )cos_k( t ) + sin_k( i )sin_k( t ),$$

then, it follows from property P4 that

$$H = \sum_{k=0}^{N-1} cos_k( i - t ),$$

and, from P9, the result follows.

## 3. The Finite Field Hartley Transform

Definition 4: Let $v = (v_0 , v_1 , ... , v_{N-1})$ be a vector of length N with components over GF(q), $q = p^r$, $p\neq 2$. The Finite Field Hartley Transform (FFHT) of v is the vector $V = (V_0 , V_1 , ... , V_{N-1})$ of components $V_k \in GI(q^m)$, given by

$$V_k := \sum_{i=0}^{N-1} v_i\, cas_k(\angle\alpha^i)$$

where $\alpha$ is a specified element of multiplicative order N in $GF(q^m)$.

Such a definition clearly mimics the classical definition of the Discrete Hartley Transform [9]. The inverse FFHT is given by the following theorem.

Theorem 2: The N-dimensional vector v can be recovered from its Hartley discrete spectrum V according to

$$v_i = \frac{1}{N(\mathrm{mod}\, p)} \sum_{k=0}^{N-1} V_k cas_k(\angle\alpha^i).$$

Proof: After substituting the $V_k$ as defined above in the expression for $v_i$, it follows that

$$v_i = \frac{1}{N(\mathrm{mod}\, p)} \sum_{k=0}^{N-1} \sum_{r=0}^{N-1} v_r cas_k(\angle\alpha^r) cas_k(\angle\alpha^i) .$$

Changing the order of summation,

$$\frac{1}{N(\bmod p)}\sum_{r=0}^{N-1} v_r \sum_{k=0}^{N-1} \mathrm{cas}_k(\angle\alpha^r)\mathrm{cas}_k(\angle\alpha^i) = \quad = \frac{1}{N(\bmod p)}\sum_{r=0}^{N-1} v_r \begin{Bmatrix} N, & i = r \\ 0, & i \neq r \end{Bmatrix} = v \qquad ð$$

A signal v and its discrete Hartley spectrum V are said to form a finite field Hartley Transform pair, denoted by $v \leftrightarrow V$. It is worthwhile to mention that the FFHT belongs to a class of discrete transforms for which the kernel of the direct and the inverse transform is exactly the same.

Letting now $g = \{g_i\} \leftrightarrow G = \{G_k\}$ and $v = \{v_i\} \leftrightarrow V = \{V_k\}$ denote FFHT pairs of length N, the following set of useful properties can be derived.

H1 - Linearity

$$ag + bv \leftrightarrow aG + bV, \quad \forall\, a,b \in GF(q).$$

H2 - Time Shift

$$\text{If } v_i = g_{i-d}\text{, then } V_k = \cos_k(d)G_k + \sin_k(d)G_{-k}.$$

H3 - Sum of Sequence (dc term)

$$V_0 = \sum_{i=0}^{N-1} v_i\,.$$

$$V_k^q = \sum_{i=0}^{N-1} v_i\ \mathrm{cas}_{N-qk}(\angle\alpha^i) = V_{N-qk}\,.$$

H4 - Initial Value

$$v_0 = \frac{1}{N(\bmod p)}\sum_{k=0}^{N-1} V_k\,.$$

H5 - Symmetry

$$G \leftrightarrow Ng\,.$$

H6 - Time Reversal

$$g_{-i} \leftrightarrow G_{-k}\,.$$

H7 - Cyclic Convolution: If ★ denotes cyclic convolution, then

$$g \star v \leftrightarrow \tfrac{1}{2}(GV + GV_{-} + G_{-}V - G_{-}V_{-}).$$

where $G_{-}$ and $V_{-}$ denotes, respectively, the sequences $\{G_{N-k}\}$ and $\{V_{N-k}\}$.

H8 - Parseval's Relation

$$N\sum_{i=0}^{N-1} g_i^2 = \sum_{k=0}^{N-1} G_k^2$$

**4. Valid Spectra**

The following lemma states a relation that must be satisfied by the components of the spectrum V for it to be a valid finite field Hartley spectrum, that is, a spectrum of a signal v with GF(q)-valued components.

Lemma 2: The vector V= $\{V_k\}$, $V_k \in GI(q^m)$, is the spectrum of a signal v = $\{v_i\}$, $v_i \in GF(q)$, $q = p^r$, if and only if

$$V_k^q = V_{N-kq}$$

where indexes are considered modulo N, i, k = 0, 1, ..., N-1 and $N \mid (q^m - 1)$.

**Proof**: From the FFHT definition and considering that $GF(p^r)$ has characteristic p, it follows that

$$V_k^q = \left(\sum_{i=0}^{N-1} v_i \operatorname{cas}_k(\angle\alpha^i)\right)^q = \left(\sum_{i=0}^{N-1} v_i^q \operatorname{cas}_k^q(\angle\alpha^i)\right).$$

If $v_i \in GF(q)\ \forall$ i, then $v_i^q = v_i$. The fact that $j^2 = -1 \notin GF(q)$ if and only if q is a prime power of the form 4s + 3, implies that $j^q = -j$. Hence,

On the other hand, suppose $V_k^q = V_{N-qk}$. Then

$$\sum_{i=0}^{N-1} v_i^q \operatorname{cas}_{N-qk}(\angle\alpha^i) = \sum_{i=0}^{N-1} v_i \operatorname{cas}_{N-qk}(\angle\alpha^i).$$

Now, let N-qk = r. Since $GCD(q^m - 1, q) = 1$, k and r ranges over the same values, which implies

$$\sum_{i=0}^{N-1} v_i^q \operatorname{cas}_r(\angle\alpha^i) = \sum_{i=0}^{N-1} v_i \operatorname{cas}_r(\angle\alpha^i),$$

r = 0, 1, ..., N-1. By the uniqueness of the FFHT, $v_i^q = v_i$ so that $v_i \in GF(q)$ and the proof is complete.

ð

Example 2 - With q = p = 3, r = 1, m = 5 and $GF(3^5)$ generated by the primitive polynomial $f(x) = x^5 + x^4 + x^2 + 1$, a FFHT of length N = 11 may be defined by taking an element or order 11 ($\alpha^{22}$ is such an element). The vectors v and V given below are an FFHT pair.

v = (1, 0, etc, etc,.........)

V = ( etc, etc, etc, ....)

The relation for valid spectra shown above implies that only two components $V_k$ are necessary to completely specify the vector V, namely $V_0$ and $V_1$. This can be verified simply by calculating the cyclotomic classes induced by lemma 1 which, in this case, are $C_0$ = (0) and $C_1$ = (1, 8, 9, 6, 4, 10, 3, 2, 5, 7).

## 5. Conclusions

In this paper, trigonometry for finite fields was introduced. In particular, the k-trigonometric functions of the *angle* of the *complex exponential* $\alpha^i$ were defined and their basic properties derived. From the $\cos_k(\angle\alpha^i)$ and $\sin_k(\angle\alpha^i)$ functions, the $\operatorname{cas}_k(\angle\alpha^i)$ function was defined and used to introduce a new Hartley Transform, the Finite Field Hartley Transform (FFHT).

The FFHT seems to have interesting applications in a number of areas. Specifically, its use in Digital Signal Processing, along the lines of the so-called number theoretic transforms (e.g. Mersenne transforms) should be investigated. In the field of error control codes, the FFHT might be used to produce a transform domain description of the field, therefore providing, possibly, an alternative to the approach introduced in [6]. Digital Multiplexing is another area that might benefit from the new Hartley Transform introduced in this paper. In particular, new schemes of efficient-bandwidth code-division-multiple-access for band-limited channels based on the FFHT are currently under development.

## Acknowledgements

The authors wish to thank Prof. James Massey for his suggestions and insightful comments which improved the final version of this paper.

## References

[1] J. M. Pollard, *The Fast Fourier Transform in a Finite Field*, Math. Comput., vol. 25, No. 114, pp. 365-374, Apr. 1971.

[2] C. M. Rader, *Discrete Convolution via Mersenne Transforms*, IEEE Trans. Comput., vol. C-21, pp. 1269-1273, Dec. 1972.

[3] I. S. Reed and T. K. Truong, *The Use of Finite Field to Compute Convolutions*, IEEE Trans. Inform. Theory, vol. IT-21, pp. 208-213, Mar. 1975.

[4] R. C. Agarwal and C. S. Burrus, *Number Theoretic Transforms to Implement Fast Digital Convolution*, IEEE Proc., vol. 63, pp. 550-560, Apr. 1975.

[5] I. S. Reed, T. K. Truong, V. S. Kwoh and E. L. Hall, *Image Processing by Transforms over a Finite Field*, IEEE Trans. Comput., vol. C-26, pp. 874-881, Sep. 1977.

[6] R. E. Blahut, *Transform Techniques for Error-Control Codes*, IBM J. Res. Dev., vol. 23, pp. 299-315, May 1979.

[7] R. M. Campello de Souza and P. G. Farrell, *Finite Field Transforms and Symmetry Groups*, Discrete Mathematics, vol. 56, pp. 111-116, 1985.

[8] J. L. Massey, *The Discrete Fourier Transform in Coding and Cryptography*, accepted for presentation at the 1998 IEEE Inform. Theory Workshop, ITW 98, San Diego, CA, Feb 9-11.

[9] R. N. Bracewell, *The Discrete Hartley Transform*, J. Opt. Soc. Amer., vol. 73, pp. 1832-1835, Dec. 1983.

[10] R. V. L. Hartley, *A More Symmetrical Fourier Analysis Applied to Transmission Problems*, Proc. IRE, vol. 30, pp. 144-150, Mar. 1942.

[11] R. N. Bracewell, *The Hartley Transform*, Oxford University Press, 1986.

[12] J.-L. Wu and J. Shiu, *Discrete Hartley Transform in Error Control Coding*, IEEE Trans. Acoust., Speech, Signal Processing, vol. ASSP-39, pp. 2356-2359, Oct. 1991.

[13] R. N. Bracewell, *Aspects of the Hartley Transform*, IEEE Proc., vol. 82, pp. 381-387, Mar. 1994.

[14] J. Hong and M. Vetterli, *Hartley Transforms Over Finite Fields*, IEEE Trans. Inform. Theory, vol. IT-39, pp. 1628-1638, Sep. 1993.